\documentclass[a4paper, 11pt]{amsart}
\usepackage[utf8]{inputenc}
\usepackage[T1]{fontenc}
\usepackage{lmodern}
\usepackage[english]{babel}
\usepackage[dvipsnames]{xcolor}
\usepackage{adjustbox}
\usepackage{letltxmacro}
\usepackage{todonotes}
\LetLtxMacro\todonotestodo\todo
\renewcommand{\todo}[2][]{\todonotestodo[#1]{TODO: {#2}}}
\usepackage{enumerate}
\usepackage{hyperref}
\usepackage{url}
\usepackage{mathrsfs}
\usepackage{amssymb}

\theoremstyle{definition}

\makeatletter
\newtheorem*{rep@theorem}{\rep@title}
\newcommand{\newreptheorem}[2]{%
\newenvironment{rep#1}[1]{%
 \def\rep@title{#2 \ref{##1}}%
 \begin{rep@theorem}}%
 {\end{rep@theorem}}}
\makeatother

\newreptheorem{theorem}{Theorem}

\newtheorem*{claim*}{Claim}

\newtheorem*{theorem*}{Theorem}
\newtheorem*{corollary*}{Corollary}
\newtheorem*{lemma*}{Lemma}

\newcommand{\Z}{\mathbb{Z}}
\newcommand{\N}{\mathbb{N}}

\title{Prime ideals in infinite products of commutative rings}
\author{Carmelo A. Finocchiaro}
\thanks{C.A. Finocchiaro is supported by GNSAGA of Istituto
Nazionale di Alta Matematica, Fondazione Cariverona (Research project
“Reducing complexity in algebra, logic, combinatorics - REDCOM” within the
framework of the program “Ricerca Scientifica di Eccellenza 2018”), and
by the Department of Mathematics and Computer Science of the University of
Catania (Research program “Proprietà algebriche locali e globali di anelli
associati a curve e ipersuperfici” PTR 2016-18).}

\author{Sophie Frisch}
\thanks{S. Frisch is supported by the Austrian Science Fund (FWF): P 35788}

\author{Daniel Windisch}
\thanks{D.~Windisch is supported by the Austrian Science Fund (FWF): P~30934}

\begin{document}

\begin{abstract}
We describe the prime ideals and, in particular, 
the maximal ideals in products $R = \prod D_\lambda$ of
families $(D_\lambda)_{\lambda \in \Lambda}$ of commutative rings. 
We show that every maximal ideal is induced by an ultrafilter 
on the Boolean algebra $\prod \mathcal{P}(\max(D_\lambda))$, 
where $\max(D_\lambda)$ is the spectrum of maximal ideals of 
$D_\lambda$, and $\mathcal{P}$ denotes the power set. 
If every $D_\lambda$ is in a certain class of rings including finite
 character domains and one-dimensional domains, we completely
characterize the maximal ideals of $R$. If every $D_\lambda$ is a 
Pr\"ufer domain, we completely characterize all prime ideals of $R$. 
\end{abstract}

\maketitle

\section{Introduction and preliminaries}

In this paper, we investigate the prime ideals of products of
commutative rings. In the case of products of Pr\"ufer rings, we
achieve a complete characterization. (Recall that Pr\"ufer rings are
the generalization of Dedekind rings defined by the property that
every finitely generated ideal is invertible).

Now prime ideals have been a mainstay of commutative ring theory ever
since Dedekind, who formalized ideas of Kummer on how to deal with
non-uniqueness of factorization of elements into irreducibles in rings
of algebraic integers, and famously found uniqueness by passing from
elements to ideals and considering factorizations of ideals into
prime ideals.

It is, therefore, an obvious priority to investigate the prime
spectrum, whenever one sets out to study a ring or a class of rings.

It is perhaps less obvious that we should study products of commutative
rings, given their awkward proclivity for having zero-divisors even when
the individual factors are integral domains.

Apart from the fact that a product is a natural object of study in
any category, we offer two reasons for studying products of rings:
first, our interest in rings of functions, that is, rings whose
elements define functions and are characterized by them, such as
rings of integer-valued polynomials~\cite{integer-valued} and rings of continuous functions~\cite{continuous}.

An identification of the elements of one ring with functions defined
on another ring clearly defines an embedding of the first ring into a
product of copies of the second ring.
If this description of ring elements as functions is a good fit,
then prime ideals of the ring of functions can be characterized by
reference to prime ideals in the product of rings. The question of
how well the interpretation of ring elements as functions is suited
to describing (finitely generated) ideals of a ring is related to
Skolem properties, which have been studied elsewhere, see~\cite{integer-valued} for an overview.

Our second reason for studying a product of rings is that it is a
step on the way to constructing ultraproducts of the same rings, and
that ultraproducts, unlike products, preserve nice properties of the
individual factors, such as, being an integral domain, a field, or,
notably, a Pr\"ufer ring. The study of the prime ideals of an infinite
product of rings is closely connected to that of the prime ideals of
all non-trivial ultraproducts of the same rings.

Some notation: Let $\Lambda$ be a set and 
$(D_\lambda)_{\lambda \in \Lambda}$ a family of commutative rings. 
Throughout this paper, we denote by 
$R = \prod_{\lambda \in \Lambda} D_\lambda$ the product of the rings 
$D_\lambda$. The prime spectrum of $\prod_{\lambda \in \Lambda} D_\lambda$ 
is closely related to a special Boolean algebra, namely:\\

\noindent
\textbf{The Boolean algebra $\prod_{\lambda \in \Lambda} \mathcal{P}(\max(D_\lambda))$,} 
which we will always denote by $\mathcal{B}$, is the product of the 
Boolean algebras $(\mathcal{P}(\max(D_\lambda)),\cap, \cup)$, 
where $\mathcal{P}(M)$ denotes the power set of a set $M$ and $\max(D)$ 
the set of all maximal ideals of a commutative ring $D$. 
Clearly, $\mathcal{B}$ is a Boolean algebra with least element 
$0_\mathcal{B} = (\emptyset)_{\lambda \in \Lambda}$. 
We denote elements $a \in R = \prod_{\lambda \in \Lambda} D_\lambda$ by 
$a = (a_\lambda) = (a_\lambda)_{\lambda \in \Lambda}$ and elements 
$Y \in \mathcal{B} =\prod_{\lambda \in \Lambda} \mathcal{P}(\max(D_\lambda))$ 
by $Y = (Y_\lambda) = (Y_\lambda)_{\lambda \in \Lambda}$.

We denote by $\dot\bigcup_{\lambda \in \Lambda} \max(D_\lambda)$ 
the disjoint union of the $\max(D_\lambda)$.
As Boolean algebras, $\mathcal{B}$ and 
$(\mathcal{P}\left(\dot\bigcup_{\lambda \in \Lambda} \max(D_\lambda)\right), \cap, \cup)$ 
are isomorphic, see Remark~\ref{3.1Remark}. 
We will make use of this fact in Section~\ref{section:primeideals}.

For an introduction to Boolean algebras, see, for instance, 
the monograph by Grätzer et al.~\cite{bool}.\\

\noindent
\textbf{Ultrafilters on Boolean algebras.}  To describe the spectrum of a 
product of commutative rings, we need ultrafilters on Boolean algebras, 
of which ultrafilters on sets are a special case.  Let $(B,\land, \lor)$ be a Boolean algebra. We denote by $0$ 
the minimal element of $B$, by $\lnot$ the complement operation on $B$ 
and by $\leq$ the canonical order relation on~$B$. For the definition of filters and ultrafilters in Boolean algebras, we again refer to the textbook by Grätzer et al.~\cite{bool}. We need these concepts in the following two special cases.
\begin{itemize}


\item[(1)] If $(B, \land, \lor) = (\mathcal{P}(\Lambda),\cap, \cup)$, 
then $0 = \emptyset$, the canonical order relation is set-theoretic 
inclusion, and $\lnot A = \Lambda \setminus A$ 
for every $A \subseteq \Lambda$. 
As is customary, we call an ultrafilter $\mathcal U$ in 
$\mathcal{P}(\Lambda)$ an \textit{ultrafilter on} $\Lambda$. 

\item[(2)] In the special case 
$B = \mathcal{B} = \prod_{\lambda \in \Lambda} \mathcal{P}(\max(D_\lambda))$,
the lattice operations on $\mathcal{B}$ are
\begin{align*}
Y \land Z = (Y_\lambda \cap \nolinebreak Z_\lambda)_{\lambda \in \Lambda} 
\quad &\text{ and } \quad Y \lor Z = (Y_\lambda \cup Z_\lambda)_{\lambda \in \Lambda},  \\
\lnot Y = (\max(D_\lambda) \setminus Y_\lambda)_{\lambda \in \Lambda} \quad &\text{ and }
 \quad 0 = 0_\mathcal{B} = (\emptyset)_{\lambda \in \Lambda}.
\end{align*}
Furthermore, $Y \leq Z$ if and only if 
$Y_\lambda \subseteq Z_\lambda$ for all $\lambda \in \Lambda$.

\item[(3)] A non-empty subset $M \subseteq B$ has the 
\textit{finite intersection property} if $Y_1 \land \hdots \land Y_n \neq 0$ 
for all $Y_1,\hdots,Y_n \in M$. If $M \subseteq B$ has the finite 
intersection property, then 
\[\mathcal{F} = \{ F \in B \mid \exists Y_1,\hdots,Y_n \in M \ Y_1 \land \hdots \land Y_n \leq F \}\] 
is easily seen to be a filter in $B$ containing $M$.

\item[(4)] Every filter in $B$ is contained in some ultrafilter in $B$. 
This follows via Zorn's Lemma from the fact that ultrafilters in $B$ are 
exactly the maximal elements with respect to set-theoretical inclusion in 
the set of all filters in $B$. 

\item[(5)] If $\mathcal{U}$ is an ultrafilter in $B$ and  $X,Y \in B$ 
such that $X \lor Y \in \mathcal{U}$, then $X \in \mathcal{U}$ or 
$Y \in \mathcal{U}$.
\end{itemize}
The above facts will be used throughout this paper 
without any additional reference.\\

In 1991, Levy, Loustaunau and Shapiro~\cite{Integers} characterized the 
maximal spectrum of an infinite product of copies of $\Z$ by ultrafilters 
on $\mathcal{B}$. They described all prime ideals of 
$R = \prod_{\lambda \in \Lambda} \Z$ and investigated the order structure 
of chains inside $\text{spec}(R)$. O'Donnell~\cite{Donnell} generalized 
some of these results to maximal ideals in products of commutative rings 
and characterized certain classes of prime ideals in products of Dedekind 
domains. Olberding, Saydam and Shapiro~\cite{Olb1,Olb3,Olb2} described 
special classes of prime ideals in ultraproducts of commutative rings 
using a totally different approach. Some of the ideas in all of these papers show similarities with approaches to the prime ideal structure of certain non-standard Dedekind domains, as treated by Cherlin~\cite{Cherlin}. This comes from the fact that such non-standard Dedekind domains are subrings of ultrapowers of Dedekind domains.

Our aim is to extend the initial idea of Levy, Loustaunau and Shapiro to 
more general situations, thereby strengthening and generalizing many of 
the above-mentioned results. To give an example, Theorem \ref{3.3} 
describes all prime ideals of arbitrary products of Prüfer domains. 
Recall that Prüfer domains are a common generalization of Dedekind 
domains and valuation domain. The usefulness of the concept is evident by 
the fact that Prüfer domains can be characterized in many different, 
natural ways. For instance, a domain is Prüfer if and only if the 
localization at every prime ideal is a valuation domain. Likewise, a 
domain is Prüfer if and only if every finitely generated ideal is 
invertible. Equivalently, a domain is Prüfer if and only if its ideal 
lattice is distributive. There are many more equivalent 
characterizations; see the monograph by Gilmer \cite[Chapter IV]{Gilmer}. 
The book by Fontana, Huckaba and Papick \cite{Pruefer} covers collected 
topics in Prüfer domains.\\

\noindent
\textbf{The Skolem-property.} A subring $T$ of $R = \prod D_\lambda$ is 
said to have the \textit{Skolem-property} if for all 
$a^{(1)},\hdots,a^{(n)} \in T$ such that the ideal 
$(a^{(1)}_\lambda,\hdots,a^{(n)}_\lambda)$ is equal to $D_\lambda$ for 
all $\lambda \in \Lambda$, it follows that 
$(a^{(1)},\hdots,a^{(n)}) = T$. \\

Note that the Skolem-property of $T$ depends not only on the 
abstract ring $T$ but also on its embedding in $R$.  \\
The Skolem-property introduced here is a generalization of the particular 
case $T = \text{Int}(D) = 
\{ f \in K[x] \mid f(D) \subseteq D \} \subseteq \prod D$, 
where $D$ is a domain with quotient field $K$. 
For a deeper insight into this topic, see 
\cite{Sk1, Sk6, Sk5, Sk3, Sk7, Sk2, Sk4}. 

In section 2, we show that $T$ has the Skolem-property if and only if 
every maximal ideal of $T$ is induced by an ultrafilter in $\mathcal{B}$. 
Note that $R = \prod D_\lambda$ itself always has the Skolem-property. 
Moreover, we characterize the ultrafilters in $\mathcal{B}$ inducing 
maximal ideals of $R$ in the case that every $D_\lambda$ satisfies the 
following property, which we call~(+):

For all $r \in D$ and $a \in D \setminus \{0\}$ there exists $d \in D$ 
such that $d$ is in every maximal ideal that contains $a$ but not $r$, 
and $d$ is in no maximal ideal that contains $r$.

We show that finite character domains and 
one-dimensional domains all satisfy (+).

We also investigate when every ultrafilter on $\mathcal{B}$ induces a 
maximal ideal of $R$. It turns out that this property is closely related 
to the proconstructibility of the maximal spectra of the component rings 
$D_\lambda$. 

Furthermore, we consider minimal prime ideals of subrings $T \subseteq R$. 
For example, we prove that every prime ideal of a product of domains $R$ 
contains exactly one minimal prime ideal. \\

\noindent
\textbf{First-order sentences and ultraproducts.} In section 3, we apply 
some model theory (in particular, \L o\'s's Theorem) to ultraproducts. We 
restrict ourselves to ultraproducts of rings and to the first-order 
language of rings whose non-logical symbols are $+$, $\cdot$, $0$ and $1$. 
Roughly speaking, a first-order sentence in this language is a formula 
without free variables only using $=$, $+$, $\cdot$, $0$, $1$, variables 
and logical symbols such as quantifiers and sentential connectives. When 
such sentences are interpreted, the underlying set is assumed to be the 
set of elements of a ring, so that quantifiers range only over ring 
elements (and not, for instance, over ideals or functions).

If $\mathcal{F}$ is an ultrafilter on $\Lambda$, we denote by 
$R^* = \prod_{\lambda \in \Lambda}^\mathcal{F} D_\lambda$ 
the ultraproduct of the $D_\lambda$, which is the quotient ring of 
$R = \prod D_\lambda$ obtained by identifying $r,s \in R$ with the 
property that $\{\lambda \in \Lambda \mid r_\lambda = s_\lambda \}$ is in 
$\mathcal{F}$.

Ultrafilters and ultraproducts are playing an
increasingly important role in commutative ring theory, for instance, in 
the work of Olberding (cf. \cite{Olb1, Olb3, Olb2}),  Fontana and Loper 
(cf. \cite{Loper2,Loper1,Loper4,Loper5,Loper3}), and Schoutens 
(cf. \cite{Schoutens2, Schoutens1}).

We will extensively use the following fundamental theorem for 
ultraproducts~\cite[Theorem 4.1.9]{mod}:

\newtheorem*{Los}{Theorem of \L o\'s}
\begin{Los}\label{Los}
A first order sentence $\varphi$ is satisfied by $R^*$ if and only if the 
set of all $\lambda \in \Lambda$ such that $D_\lambda$ satisfies $\varphi$ 
is in $\mathcal{F}$.
\end{Los}

From the Theorem of \L o\'s, it follows in particular that, if every 
$D_\lambda$ is an integral domain (respectively a field), then so is $R^*$. 
Moreover, it can be easily seen that the quotient field of $R^*$ is 
isomorphic to the ultraproduct of the quotient fields of the $D_\lambda$. 
For a more precise and general treatment of the introduced concepts, 
see~\cite{mod}.\\

It is known that being a Prüfer domain is preserved by 
ultraproducts~\cite[Proposition 2.2]{Olb1}. 
In section 3, we explicitly describe the valuation on the quotient field 
$K^*$ of an ultraproduct $R^*$ of Prüfer domains $D_\lambda$ whose 
valuation ring is the localization $R^*_M$ at a maximal ideal 
$M \subseteq R^*$.

Also, in the case when  all $D_\lambda$ are Prüfer domains, we are able 
to describe all prime ideals in $R= \prod D_\lambda$, using a common 
generalization of concepts introduced in~\cite{Integers} and~\cite{Olb3}. 
It turns out that every non-minimal prime ideal of $R$ contained in a 
maximal ideal of a certain type is of infinite height. We show that 
maximal ideals of this type always exist. Finally, we construct examples 
of Pr\"ufer domains such that every non-zero prime ideal is of infinite 
height. The Prüfer domain can be chosen either local 
(and therefore a valuation domain) or non-local.

\section{Maximal ideals and minimal prime ideals}\label{section:maxideals}

\subsection*{Describing all maximal ideals}

Let $D$ be a commutative ring. We denote the spectrum of all maximal ideals of $D$ by $\max(D)$. For an ideal $I \subseteq D$, we write $\mathscr{V}(I)$ for the set of all maximal ideals of $D$ containing $I$ and $\mathscr{D}(I)$ for $\max(D) \setminus \mathscr{V}(I)$. If $I = (a_1,\hdots,a_n)$ is finitely generated, we let $ \mathscr V(a_1,\hdots,a_n) = \mathscr{V}(I)$.\\
For an element $a \in R = \prod D_\lambda$, we set $S(a) = (\mathscr V(a_\lambda))_{\lambda \in \Lambda}$, an element of $\mathcal{B} = \prod \mathcal{P}(\max(D_\lambda))$. Moreover, if $\mathcal{U}$ is a filter in $\mathcal{B}$ and $T $ a subring of $R$, we define
\begin{align*}
(\mathcal{U})_T = \{ a \in T \mid S(a) \in \mathcal{U} \}.
\end{align*}
We omit $T$ whenever it is understood from the context.

There is one trivial case, where this set $(\mathcal{U})_T$ can be understood easily: Pick some $\mu \in \Lambda$ and some maximal ideal $M$ of $D_\lambda$. Then the element $S = (S_\lambda)_{\lambda \in \Lambda}$, where $S_\mu = \{M\}$ and $S_\lambda = \emptyset$ for $\lambda \neq \mu$, defines the ultrafilter $\mathcal{U}_S = \{B \in \mathcal{B} \mid S \leq B\}$, the \textit{principal ultrafilter} generated by $S$.

With this notation, $(\mathcal{U}_S)_R$ equals $\prod_{\lambda \in \Lambda} P_\lambda$, where $P_\mu = M$ and $P_\lambda = R_\lambda$ for $\lambda \neq \mu$, and this set is clearly a prime ideal of $R$. The next lemma shows that this last statement is true in full generality.

\newtheorem{2.1}{Lemma}[section]
\begin{2.1}\label{2.1}
Let $T \subseteq R$ be a subring, $a,b \in R$ and $\mathcal{U}$ be a filter in $\mathcal{B}$. Then the following assertions hold:
\begin{itemize}
\item[(1)] $S(a) \land S(b) = (\mathscr V(a_\lambda,b_\lambda))_{\lambda \in \Lambda}$.
\item[(2)] $S(a) \lor S(b) = S(ab)$.
\item[(3)] $(\mathcal{U})_T$ is an ideal of $T$.
\item[(4)] If $\mathcal{U}$ is an ultrafilter in $\mathcal{B}$, then $(\mathcal{U})_T$ is a prime ideal of $T$.
\end{itemize}
\end{2.1}

\begin{proof}
(1), (2) and (3) follow immediately from the relevant definitions. \\
For the proof of (4), let $\mathcal{U}$ be an ultrafilter in $\mathcal{B}$ and note that $ 1 \notin (\mathcal{U})$, because $S(1) = 0_\mathcal{B} \notin \mathcal{U}$. If now $a,b \in T$ such that $ab \in (\mathcal{U})$, then by (2) we have that $S(a) \lor S(b) = S(ab) \in \mathcal{U}$. Since $\mathcal{U}$ is an ultrafilter, it follows that $S(a) \in \mathcal{U}$ or $S(b) \in \mathcal{U}$ and therefore $a \in (\mathcal{U})$ or $b \in (\mathcal{U})$.
\end{proof}

\newtheorem{2.2}[2.1]{Proposition}
\begin{2.2}\label{2.2}
For a subring $T \subseteq R = \prod_{\lambda \in \Lambda} D_\lambda$ the following assertions are equivalent:
\begin{itemize}
\item[(a)] $T$ has the Skolem-property.
\item[(b)] For every proper ideal $\mathfrak{A} \subseteq T$ the set $\{S(a) \mid a \in \mathfrak{A}\} \subseteq \mathcal{B}$ satisfies the finite intersection property.
\item[(c)] Every proper ideal of $T$ is contained in an ideal of the form $(\mathcal{U})$, where $\mathcal{U}$ is an ultrafilter in $\mathcal{B}$.
\item[(d)] Every maximal ideal of $T$ is of the form $(\mathcal{U})$ for some ultrafilter $\mathcal{U}$ in $\mathcal{B}$.
\end{itemize}
\end{2.2}

\begin{proof}
"(a) $\Rightarrow$ (b)": Let $\mathfrak{A} \subseteq T$ be a proper ideal and $a^{(1)},\hdots,a^{(n)} \in \mathfrak{A}$. Then, by Lemma \ref{2.1}(1), it follows that $S(a^{(1)}) \land \hdots\land S(a^{(n)}) = (\mathscr V(a^{(1)}_\lambda,\hdots,a^{(n)}_\lambda))$. Assume to the contrary that $\mathscr V(a^{(1)}_\lambda,\hdots,a^{(n)}_\lambda) = \emptyset$ for all $\lambda \in \Lambda$. Then $(a^{(1)}_\lambda,\hdots,a^{(n)}_\lambda) = D_\lambda$ for all $\lambda \in \Lambda$ and by the Skolem-property $\mathfrak{A} \supseteq (a^{(1)},\hdots,a^{(n)}) = T$, which is a contradiction. \\
"(b) $\Rightarrow$ (c)": Let $\mathfrak{A} \subseteq T$ be a proper ideal. Then, by (b), we can pick an ultrafilter $\mathcal{U}$ in $\mathcal{B}$ such that $\{S(a) \mid a \in \mathfrak{A}\} \subseteq \mathcal{U}$. Now it follows by definition that $\mathfrak{A} \subseteq (\mathcal{U})$. \\
"(c) $\Rightarrow$ (d)": This is clear. \\
"(d) $\Rightarrow$ (a)": Let $a^{(1)},\hdots,a^{(n)} \in T$ such that $\mathfrak{A} = (a^{(1)},\hdots,a^{(n)})$ is a proper ideal of $T$. Let $\mathcal{U}$ be an ultrafilter in $\mathcal{B}$ such that $\mathfrak{A} \subseteq (\mathcal{U})$. We want to show that $(a^{(1)}_\lambda,\hdots,a^{(n)}_\lambda)$ is proper for some $\lambda \in \Lambda$. Assume to contrary that $(a^{(1)}_\lambda,\hdots,a^{(n)}_\lambda) = D_\lambda$ for all $\lambda \in \Lambda$. Then $0_\mathcal{B} = (\mathscr V(a^{(1)}_\lambda,\hdots,a^{(n)}_\lambda)) = S(a^{(1)}) \land \hdots\land S(a^{(n)}) \in \mathcal{U}$, which is a contradiction.
\end{proof}

\newtheorem{2.3}[2.1]{Corollary}
\begin{2.3}\label{2.3}
Let $(D_\lambda)_{\lambda \in \Lambda}$ be a family of commutative rings. Then every maximal ideal of $R = \prod_{\lambda \in \Lambda} D_\lambda$ is of the form $(\mathcal{U})$ for some ultrafilter $\mathcal{U}$ in $\mathcal{B} = \prod_{\lambda \in \Lambda} \mathcal{P}(\max(D_\lambda))$.
\end{2.3}

\subsection*{Characterizing ultrafilters that induce maximal ideals}

\theoremstyle{definition}
\newtheorem{2.4}[2.1]{Definition}
\begin{2.4}\label{2.4}
We say that a ring $D$ satisfies \textit{property} (+) if for all $r \in D$ and $a \in D \setminus \{0\}$ there exists $d \in D$ such that $d$ is in every maximal ideal that contains $a$ but not $r$, and $d$ is in no maximal ideal that contains $r$.
\end{2.4}

We will see that property (+) gives a setting where we can characterize those ultrafilters in $\mathcal{B}$ that induce maximal ideals of $R$.

We first give some easy equivalences to property (+), which will help us to give examples of classes of domains satisfying or not satisfying it. We need the following fact.

\theoremstyle{definition}
\newtheorem{2.4a}[2.1]{Remark}
\begin{2.4a}\label{2.4a}

 If $I$ is an ideal of a ring $D$ and $r \in D$, then $I \subseteq \bigcup_{Q \in \mathscr V(r)} Q$ implies that $I \subseteq Q$ for some $ Q \in \mathscr V(r)$. Indeed, $I \subseteq \bigcup_{Q \in \mathscr V(r)} Q$ is equivalent to $iD+rD \neq D$ for every $i \in I$, and $I \subseteq Q$ for some $Q \in \mathscr V(r)$ is equivalent to $I + rD \neq D$, which is the same thing. This is an instance of (potentially) infinite prime avoidance, a phenomenon investigated in more detail by J.~Chen~\cite{Chen}.
 
 \end{2.4a}

\newtheorem{2.5}[2.1]{Lemma}
\begin{2.5}\label{2.5}
Let $D$ be a ring. Then the following are equivalent:
\begin{itemize}
\item[(a)] $D$ satisfies property (+).
\item[(b)] For all $r \in D$, $a \in D \setminus \{0\}$ and $Q \in \mathscr V(r)$ we have that $(\bigcap_{M \in \mathscr V(a) \setminus \mathscr V(r)} M) \setminus Q \neq \emptyset$.
\item[(c)] For all $r \in D$, $a \in D \setminus \{0\}$ and all maximal ideals $Q \subseteq D$ such that $\bigcap_{M \in \mathscr V(a) \setminus \mathscr V(r)} M \subseteq Q$ it follows that there exists some $M \in \mathscr V(a) \setminus \mathscr V(r)$ such that $M \subseteq Q$.
\item[(d)]  For all $r \in D$, $a \in D \setminus \{0\}$ and all maximal ideals $Q \subseteq D$ such that $\bigcap_{M \in \mathscr V(a) \setminus \mathscr V(r)} M \subseteq Q$ it follows that $Q \in \mathscr D(r)$.
\item[(e)] For all $r \in D$ and for all non-zero $a \in D$, there exists $d \in D$ such that \[\mathscr V(a) \cap \mathscr D(r) \subseteq \mathscr V(d) \subseteq \mathscr D(r)\].
\end{itemize}

\end{2.5}

\begin{proof}
The equivalence of (b), (c) and (d) is clear. Also (b) follows immediately from (a). Moreover, (a) and (e) are trivially equivalent. It now suffices to prove "(b) $\Rightarrow$ (a)". Assume that (a) does not hold. Then $\bigcap_{M \in \mathscr V(a) \setminus \mathscr V(r)} M \subseteq  \bigcup_{Q \in \mathscr V(r)} Q$ for some $r \in D$ and some non-zero $a \in D$. By Remark~\ref{2.4a}, it follows that $\bigcap_{M \in \mathscr V(a) \setminus \mathscr V(r)} M \subseteq Q$ for some $Q \in \mathscr V(r)$, which contradicts (b).
\end{proof}

\newtheorem{2.6}[2.1]{Example}
\begin{2.6}\label{2.6}
Let $D$ be a domain of finite character, i.e., every $a \in D \setminus \{0\}$ is contained in only finitely many maximal ideals of $D$. It is immediate by (c) in Lemma \ref{2.5} and the fact that every maximal ideal $Q \subseteq D$ is prime that $D$ satisfies (+).\\
In particular, principal ideal domains and, more generally, one-dimensional Noetherian domains satisfy (+).
\end{2.6}

In the case that $D$ does not have finite character, the situation is much more involved, as we want to illustrate by the next example. Nevertheless, Proposition \ref{2.8} will enlarge the class of rings of which we know that they satisfy (+) into an important direction.

\newtheorem{2.7}[2.1]{Example}
\begin{2.7}\label{2.7}
\begin{itemize}
\item[(1)] If $K$ is a field and $n \geq 2$, then the polynomial ring in $n$ indeterminates over $K$ is a Noetherian factorial domain of Krull dimension $n$ that is not Prüfer and does not satisfy property (+).
\item[(2)] The polynomial ring $\mathbb{Z}[x]$ is a two-dimensional Noetherian factorial domain that is not Prüfer and does not satisfy property (+).
\item[(3)] The ring of integer-valued polynomials $\text{Int}(\mathbb{Z})$ is a two-dimensional non-Noetherian Prüfer domain not satisfying (+).
\end{itemize}

\end{2.7}

\newtheorem{2.8}[2.1]{Proposition}
\begin{2.8}\label{2.8}
Every one-dimensional domain satisfies (+).
\end{2.8}

\begin{proof}
Let $D$ be one-dimensional, $r \in D$ and $a \in D \setminus \{0\}$. Note that $Z = \bigcap_{M \in \mathscr V(a) \cap \mathscr D(r)} M$ is an intersection of prime ideals with $a \in Z$. Therefore $Z$ is a non-zero radical ideal of $D$. If $Z = D$, the assertion is trivial, so assume that $Z$ is a proper ideal, which implies that $D/Z$ is a reduced zero-dimensional ring (i.e., von Neumann regular).

Since $(r+Z)$ is a principal ideal of $D/Z$, there exists $e \in D$ such that $e + Z$ is idempotent in $D/Z$ and $(r+Z) = (e+Z)$. We define $d = 1 - e$ and claim that $d$ is the right choice for property~(+). 

Let $M \subseteq D$ be a maximal ideal such that $ a \in M$ and $r \notin M$. This implies $r+Z \notin M/Z$. (If we had $r + Z \in M/Z$, then we could pick $m \in M$ such that $r + Z = m + Z$. But then $r - m \in Z \subseteq M$, which would imply $r \in M$, a contradiction.) It follows that $e + Z \notin M/Z$ and therefore $d + Z \in M/Z$. With the same argument as before, we get $d \in M$.

Now let $Q \subseteq D$ be a maximal ideal containing $d$. Then $d + Z \in Q/Z$. Therefore $e+Z \notin Q/Z$, which implies $r+Z \notin Q/Z$ and hence $r \notin Q$.
\end{proof}

Note that Proposition \ref{2.8} gives also rise to examples of domains satisfying (+) that are not of finite character. For instance, let $\bar{\mathbb{Z}}$ be the integral closure of $\mathbb{Z}$ in some algebraic closure of $\mathbb{Q}$. Then $\bar{\mathbb{Z}}$ is a one-dimensional Prüfer domain but it is not of finite character. Indeed, every prime number $p \in \mathbb{Z}$ is contained in infinitely many maximal ideals of $\bar{\mathbb{Z}}$.

We now return to the investigation of maximal ideals of $R = \prod D_\lambda$ and ultrafilters in $\mathcal{B} = \prod \mathcal{P}(\max(D_\lambda))$.

\newtheorem{2.9}[2.1]{Proposition}
\begin{2.9}\label{2.9}
Let $(D_\lambda)_{\lambda \in \Lambda}$ be a family of rings satisfying (+) and $\mathcal{U}$ an ultrafilter in $\mathcal{B}$ containing an element of the form $(\mathscr V(a_\lambda))_{\lambda \in \Lambda}$, where $a_\lambda \in D_\lambda \setminus \{0\}$ for all $\lambda \in \Lambda$. Then $(\mathcal{U})$ is a maximal ideal of $R = \prod D_\lambda$.
\end{2.9}

\begin{proof}
Let $r \in R \setminus (\mathcal{U})$ and let $(a_\lambda)_{\lambda \in \Lambda}$ be a family such that $a_\lambda \in D_\lambda \setminus \{0\}$ for all $\lambda \in \Lambda$ and $(\mathscr V(a_\lambda))_{\lambda \in \Lambda} \in \mathcal{U}$. Since every $D_\lambda$ satisfies (+), for each $\lambda \in \Lambda$ we can pick some $d_\lambda \in D_\lambda$ such that $\mathscr V(a_\lambda) \cap \mathscr D(r_\lambda) \subseteq \mathscr V(d_\lambda) \subseteq \mathscr D(r_\lambda)$ and define $d = (d_\lambda)_{\lambda \in \Lambda}$. Since $r \notin (\mathcal{U})$, it follows that $S(r) \notin \mathcal{U}$ and therefore $(\mathscr D(r_\lambda)) = \lnot S(r) \in \mathcal{U}$, because $\mathcal{U}$ is an ultrafilter. Hence we have $S(d) \geq S(a) \land (\mathscr D(r_\lambda)) \in \mathcal{U}$, which implies that $S(d) \in \mathcal{U}$ and therefore $d \in (\mathcal{U})$. On the other hand, we have $(d_\lambda,r_\lambda) = D_\lambda$ for all $\lambda \in \Lambda$. By the Skolem-property of $R$ it follows that $(d,r) = R$ and therefore $(\mathcal{U})$ is maximal.
\end{proof}

We now introduce two new kinds of ideals. The first one will also be the prototype of minimal prime ideals in subrings $T \subseteq R$.

\newtheorem{2.9a}[2.1]{Definition}
\begin{2.9a}\label{2.9a}
Let $\mathcal{F}$ be an ultrafilter on $\Lambda$ and $T \subseteq R$ be a subring. For an element $x \in T$ we set $z(x) = \{\lambda\in \Lambda \mid x_\lambda = 0 \}$ and we define
\begin{align*}
(0)_\mathcal{F}^T = \{x \in T \mid z(x) \in \mathcal{F}\}.
\end{align*}
Moreover, for a family $M = (M_\lambda)_{\lambda \in \Lambda}$, where $M_\lambda \in \max(D_\lambda)$ for every $\lambda \in \Lambda$, we set $z_M(x) = \{\lambda \in \Lambda \mid x_\lambda \in M_\lambda \}$ for an element $x \in T$ and define
\begin{align*}
M_\mathcal{F}^T = \{x \in T \mid z_M(x) \in \mathcal{F}\}.
\end{align*}
We write $(0)_\mathcal{F}^T =(0)_\mathcal{F}$ and $M_\mathcal{F}^T = M_\mathcal{F}$ whenever $T$ is understood from the context.
\end{2.9a}

Note that $(0)_\mathcal{F}^T$ and $M_\mathcal{F}^T$ are ideals of $T$.

\newtheorem{2.10}[2.1]{Lemma}
\begin{2.10}\label{2.10}
Let $\mathcal{F}$ be an ultrafilter on $\Lambda$ and $T \subseteq R$ a subring such that there exists $c \in T$ where $c_\lambda \in D_\lambda$ is a non-zero non-unit for every $\lambda \in \Lambda$. Then $(0)_\mathcal{F}$ is a non-maximal ideal of $T$. 
\end{2.10}

\begin{proof}
Let $c \in T$ as in the assumption of the lemma and let $M = (M_\lambda)$ be a family such that each $M_\lambda$ is a maximal ideal of $D_\lambda$ containing $c_\lambda$. Clearly, $M_\mathcal{F} \subseteq T$ is a proper ideal with $(0)_\mathcal{F} \subseteq M_\mathcal{F}$ and $c \in M_\mathcal{F} \setminus (0)_\mathcal{F}$. Therefore $(0)_\mathcal{F}$ is not maximal.
\end{proof}

\newtheorem{2.11}[2.1]{Proposition}
\begin{2.11}\label{2.11}
Let $R = \prod_{\lambda \in \Lambda} D_\lambda$, where each $D_\lambda$ is a domain, and $T$ a subring of $R$ with the property that there exists $c \in T$ such that $c_\lambda \in D_\lambda$ is a non-zero non-unit for every $\lambda \in \Lambda$. Moreover, let $\mathcal{U}$ be an ultrafilter in $\mathcal{B}$ such that $(\mathcal{U}) \subseteq T$ is a maximal ideal. Then $\mathcal{U}$ contains an element of the form $(\mathscr V(a_\lambda))_{\lambda \in \Lambda}$, where $a_\lambda \in D_\lambda \setminus \{0\}$ for all $\lambda \in \Lambda$.
\end{2.11}

\begin{proof}
First, note that $\{z(x) \mid x \in (\mathcal{U})\}$ does not have the finite intersection property. For, otherwise there would exist an ultrafilter $\mathcal{F}$ on $\Lambda$ such that $(\mathcal{U}) \subseteq (0)_\mathcal{F}$, which would imply that $(0)_\mathcal{F}$ is maximal, contradiction to Lemma \ref{2.10}.

So we can pick $x^{(1)},\hdots,x^{(n)} \in (\mathcal{U})$ such that $z(x^{(1)}) \cap \hdots \cap z(x^{(n)}) = \emptyset$. Therefore for all $\lambda \in \Lambda$ we can choose $i_\lambda \in \{1,\hdots,n\}$ such that $x_\lambda^{(i_\lambda)} \neq 0$ and therefore $a_\lambda := c_\lambda \cdot x_\lambda^{(i_\lambda)}$ is a non-zero non-unit of $D_\lambda$. If we now set $a = (a_\lambda)_{\lambda \in \Lambda}$, then $(\mathscr V(a_\lambda)) = (\mathscr V(c_\lambda \cdot x_\lambda^{(i_\lambda)})) \geq (\mathscr V(c_\lambda \cdot x^{(1)}_\lambda,\hdots,c_\lambda \cdot x^{(n)}_\lambda)) = S(c \cdot x^{(1)}) \land \hdots \land S(c \cdot x^{(n)}) \in \mathcal{U}$. Therefore $(\mathscr V(a_\lambda)) \in \mathcal{U}$, which we wanted to show.
\end{proof}

\newtheorem{2.12}[2.1]{Corollary}
\begin{2.12}\label{2.12}
Let $(D_\lambda)_{\lambda \in \Lambda}$ be a family of domains satisfying property (+) none of which is a field and $R = \prod_{\lambda \in \Lambda} D_\lambda$. Then the maximal ideals of $R$ are exactly the ideals of the form $(\mathcal{U})$, where $\mathcal{U}$ is an ultrafilter in the Boolean algebra $\mathcal{B} = \prod_{\lambda \in \Lambda} \mathcal{P}(\max(D_\lambda))$ containing an element of the form $(\mathscr V(a_\lambda))_{\lambda \in \Lambda}$ such that $a_\lambda \neq 0$ for all $\lambda \in \Lambda$.
\end{2.12}

\subsection*{The finite character case}

Recall that a ring $D$ is said to be of finite character if every non-zero element of $D$ is contained in at most finitely many maximal ideals.
Let $(D_\lambda)_{\lambda \in \Lambda}$ be a family of domains of finite character and $T \subseteq R = \prod_{\lambda \in \Lambda} D_\lambda$ a subring containing some $c = (c_\lambda)_{\lambda \in \Lambda}$ such that every $c_\lambda \in D_\lambda$ is a non-zero non-unit. It follows immediately from Proposition \ref{2.11} that every ultrafilter $\mathcal{U}$ such that $(\mathcal U) \subseteq T$ is a maximal ideal must contain an element $Y = (Y_\lambda)$ where each $Y_\lambda$ is finite.

The next result yields a statement analogous to Proposition \ref{2.9} in the finite character case.

\theoremstyle{definition}
\newtheorem{2.13}[2.1]{Proposition}
\begin{2.13}\label{2.13}
Let $(D_\lambda)_{\lambda \in \Lambda}$ be a family of rings such that for every $\lambda \in \Lambda$ the ring $D_\lambda/J(D_\lambda)$ is of finite character, where $J(D_\lambda)$ stands for the Jacobson radical of $D_\lambda$. Let $\mathcal{U}$ be an ultrafilter in $\mathcal{B}$ containing an element $Y = (Y_\lambda)_{\lambda \in \Lambda}$ such that $Y_\lambda$ is finite for every $\lambda \in \Lambda$. Then $(\mathcal{U})$ is a maximal ideal of $R = \prod_{\lambda \in \Lambda} D_\lambda$.
\end{2.13}

\begin{proof}
Let $r \in R \setminus (\mathcal{U})$. We show $(r) + (\mathcal U) = R$. Choose $a = (a_\lambda)$ such that 
\begin{itemize}
\item[(1)] $a_\lambda \in P$ for all $P \in \mathscr D(r_\lambda) \cap Y_\lambda$ and
\item[(2)] $a_\lambda \notin Q$ for all $Q \in \mathscr V(r_\lambda)$.
\end{itemize}
If $\mathscr V(r_\lambda)$ is finite then this is possible by the Chinese Remainder Theorem. If $\mathscr V(r_\lambda) = \max(D_\lambda)$ then setting $a_\lambda = 1$ works. To see that $a \in (\mathcal{U})$, note that $S(r) \notin \mathcal{U}$ and therefore $(\mathscr D(r_\lambda)) = \lnot S(r) \in \mathcal{U}$. It follows that $S(a) \geq (\mathscr D(r_\lambda)) \land Y \in \mathcal{U}$ and hence $a \in (\mathcal{U})$. By (2) and the Skolem-property of $R$, it follows that $(a,r) = R$.
\end{proof}

\newtheorem{2.14}[2.1]{Corollary}
\begin{2.14}\label{2.14}
Let $(D_\lambda)_{\lambda \in \Lambda}$ be a family of domains of finite character none of which is a field and let $R = \prod_{\lambda \in \Lambda} D_\lambda$.
Then the maximal ideals of $R$ are exactly the ideals of the form $(\mathcal{U})$, where $\mathcal{U}$ is an ultrafilter in the Boolean algebra $B = \prod_{\lambda \in \Lambda} \mathcal{P}(\max(D_\lambda))$ containing an element $Y = (Y_\lambda)_{\lambda \in \Lambda}$ such that $Y_\lambda$ is finite for all $\lambda \in \Lambda$.
\end{2.14}

\subsection*{Proconstructibility of the maximal spectra}

We now investigate the connection between a topological property of $\max(D_\lambda)$ called proconstructibility and the property of $R = \prod D_\lambda$ that $(\mathcal U) \subseteq R$ is a maximal ideal for every ultrafilter $\mathcal{U}$ in $\mathcal{B}$.

If $D$ is a commutative ring, then the constructible topology on $\text{spec}(D)$ is the coarsest topology for which all open and compact subsets (with respect to the Zariski topology) are clopen. This topology makes $\text{spec}(D)$ a compact Haussdorf space and preserves important properties. The closed sets of constructible topology are called \textit{proconstructible} sets. To describe proconstructible sets, one can use the fact that the constructible topology is equal to the so-called ultrafilter topology on $\text{spec}(D)$: A subset $X \subseteq \text{spec}(D)$ is proconstructible if and only if for each ultrafilter $\mathcal U$ on $X$ the prime ideal $X_\mathcal{U} = \{ r \in D \mid V(r) \cap X \in \mathcal{U} \}$ of $D$ is in $X$, where $V(r) = \{ P \in \text{spec}(D) \mid r \in P \}$. If we consider the subspace $X = \max(D)$ of $\text{spec}(D)$, then this property translates as follows: $X = \max(D)$ is proconstructible if and only if $X_\mathcal{U} = \{ r \in D \mid \mathscr V(r) \in \mathcal{U} \}$ is maximal for each ultrafilter $\mathcal{U}$ on $\max(D)$.

Note also that a subset $Y$ of $\text{spec}(D)$ is proconstructible if and only if $Y$ is a compact subset with respect to the constructible topology. This is because $\text{spec}(D)$ with the constructible topology is a compact Haussdorf space.

It is also worth mentioning that for a domain $D$ of finite character that has infinitely many maximal ideals, $\max(D)$ is never proconstructible in $\text{spec}(D)$. Indeed the closure of $\max(D)$ with respect to the constructible topology is $\max(D) \cup \{(0)\}$, see the paper of the first author and Tartarone~\cite[Lemma 2.9]{Carmelo-Tartarone}.

\theoremstyle{definition}
\newtheorem{2.15}[2.1]{Proposition}
\begin{2.15}\label{2.15}
If $(\mathcal{U})$ is a maximal ideal of $R = \prod_{\lambda \in \Lambda} D_\lambda$ for every ultrafilter $\mathcal{U}$ in $\mathcal{B}$, then $\max(D_\lambda)$ is proconstructible in $\text{spec}(D_\lambda)$ for every $\lambda \in \Lambda$.
\end{2.15}

\begin{proof}
Fix $\lambda \in \Lambda$ and set $X = \max(D_\lambda)$. As noted before the proposition, it suffices to show that $X_F = \{ r \in D_\lambda \mid \mathscr V(r) \in F \}$ is in $X$ for every ultrafilter $F$ on $X$. So let $F$ be an ultrafilter on $X$. For every $r \in D_\lambda$ consider the element $Y^{(r)} \in \mathcal{B}$ defined by setting
\begin{align*}
Y^{(r)}_\mu &= \mathscr D(r) \text{ if } \mu = \lambda\\
Y^{(r)}_\mu &= \emptyset \text{ if } \mu \neq \lambda
\end{align*}
for $\mu \in \Lambda$.

Now consider the subset $\mathcal{G} = \{Y^{(r)} \mid r \in D_\lambda \setminus X_F \}$ of $\mathcal{B}$. Since $F$ is an ultrafilter on $X = \max(D_\lambda)$ and $\mathscr V(r) \notin F$ for every $r \in D_\lambda \setminus X_F$, it follows that $\mathscr D(r_1) \cap \hdots \cap \mathscr D(r_n) \in F$ for all $r_1,\hdots,r_n \in D_\lambda \setminus X_F$. Hence $\mathcal{G}$ has the finite intersection property as a subset of the Boolean algebra $\mathcal{B}$. Let $\mathcal{U}$ be an ultrafilter in $\mathcal{B}$ such that $\mathcal{G} \subseteq \mathcal{U}$.

By assumption $(\mathcal{U}) \subseteq R$ is a maximal ideal and it can easily be seen that it contains the kernel of the projection map $p: R \to D_\lambda$. Indeed, if $r \in R$ such that $r_\lambda = 0$, then $S(r) \geq Y^{(1)} \in \mathcal{U}$. It follows that $p((\mathcal{U})) \subseteq D_\lambda$ is a maximal ideal.

Now we claim that $D_\lambda \setminus X_F \subseteq D_\lambda \setminus p((\mathcal{U}))$. If we know this, it follows that $p((\mathcal{U})) \subseteq X_F$ and therefore $X_F$ is maximal, which is what we wanted to show.

To prove the claim, assume to the contrary that there exists $\alpha \in D_\lambda \setminus X_F$ such that $\alpha = p(f)$ for some $f \in (\mathcal{U})$, i.e., $\alpha = f_\lambda$. Since $S(f)$ and $Y^{(\alpha)}$ are in $\mathcal{U}$, it follows that $0_\mathcal{B} = S(f) \land Y^{(\alpha)} \in \mathcal{U}$, which is a contradiction.
\end{proof}

\newtheorem{2.16}[2.1]{Definition}
\begin{2.16}\label{2.16}
A commutative ring $D$ is said to satisfy \textit{property} (++), if for all $r \in D$ there exists some $d \in D$ such that $\mathscr D(r) = \mathscr V(d)$.
\end{2.16}

Note that if a ring $D$ satisfies (++), then it also satisfies (+). Indeed, given $r \in D$ and $a \in D \setminus \{0\}$, let $d \in D$ such that $\mathscr D(r) = \mathscr V(d)$. Then $\mathscr V(a) \cap \mathscr D(r) \subseteq \mathscr D(r) = \mathscr V(d) \subseteq \mathscr D(r)$. So $D$ satisfies (+) by Lemma \ref{2.5}.

Before giving examples of rings with property (++), we show that it is relevant to our description of maximal ideals.

\newtheorem{2.17}[2.1]{Lemma}
\begin{2.17}\label{2.17}
If $(D_\lambda)_{\lambda \in \Lambda}$ is a family of commutative rings satisfying (++) then $(\mathcal{U})$ is a maximal ideal of $R = \prod_{\lambda \in \Lambda} D_\lambda$ for every ultrafilter $\mathcal{U}$ in $\mathcal{B}$.
\end{2.17}

\begin{proof}
Let $\mathcal{U}$ be an ultrafilter in $\mathcal{B}$ and choose $r \in R \setminus (\mathcal{U})$. Using property (++), let $d \in R$ such that $\mathscr D(r_\lambda) = \mathscr V(d_\lambda)$ for every $\lambda \in \Lambda$. Then $d \in (\mathcal{U})$ since $S(r) \notin \mathcal{U}$ implies $S(d) = (\mathscr V(d_\lambda)) = (\mathscr D(r_\lambda)) = \lnot S(r) \in \mathcal{U}$.  Also, $(r,d) = R$ by the Skolem-property of $R$. This shows that $(\mathcal{U})$ is maximal.
\end{proof}

\newtheorem{2.17a}[2.1]{Definition}
\begin{2.17a}\label{2.17a}
For a subset $X \subseteq \text{spec}(D)$, where $D$ is a commutative ring, we denote by $Cl^{zar}(X)$ the closure of $X$ with respect to the Zariski topology, by $Cl^{cons}(X)$ the closure of $X$ with respect to the constructible topology and by 
\begin{align*}
X^{sp} = \{ P \in \text{spec}(D) \mid P \supseteq Q \text{ for some } Q \in X \}
\end{align*}
the \textit{specialization of} $X$.
\end{2.17a}

\newtheorem{2.17b}[2.1]{Remark}
\begin{2.17b}\label{2.17b}
Fontana~\cite[Lemma 1.1]{Fontana} showed that $Cl^{zar}(X) = (Cl^{cons}(X))^{sp}$ for every $X \subseteq \text{spec}(D)$. From this it follows easily that $\max(D)$ is proconstructible in $\text{spec}(D)$ if and only if it is closed with respect to the Zariski topology on $\text{spec}(D)$. 
\end{2.17b}

\newtheorem{2.18}[2.1]{Proposition}
\begin{2.18}\label{2.18}
Let $D$ be a commutative ring such that $\max(D)$ is proconstructible in $\text{spec}(D)$. Then $D$ satisfies property (++).
\end{2.18}

\begin{proof}
Let $J$ denote the Jacobson radical of $D$. Since $\max(D)$ is proconstructible, it follows by Remark~\ref{2.17b} that $\max(D)$ is closed with respect to the Zariski topology. In this case we have that $\{P \in \text{spec}(D) \mid J \subseteq P \} = \max(D)$ and therefore $D' := D/J$ is a zero-dimensional reduced ring.

Let $r \in D$. Since $D'$ is zero-dimensional reduced, there exists some $e \in D$ such that $e + J \in D'$ is idempotent and the principal ideals $(r+J)D'$ and $(e+J)D'$ coincide. Let $d := 1 - e$. Then it can be easily seen that $\mathscr D(r+J) = \mathscr D(e+J) = \mathscr V(d+J)$. From this it is clear that $\mathscr D(r) = \mathscr V(d)$.
\end{proof}

Note that, if $D$ is zero-dimensional, then $\max(D) = \text{spec}(D)$ is proconstructible. Also, if $D$ is a one-dimensional domain with non-zero Jacobson radical $J$, then $\max(D) = V(J)$ is proconstructible. Hence both zero-dimensional rings and one-dimensional domains with non-zero Jacobson radical satisfy (++).

The next result is now an immediate consequence of Proposition \ref{2.15}, Lemma \ref{2.17} and Proposition \ref{2.18}.

\newtheorem{2.19}[2.1]{Corollary}
\begin{2.19}\label{2.19}
Let $(D_\lambda)_{\lambda \in \Lambda}$ be a family of commutative rings and $R = \prod_{\lambda \in \Lambda} D_\lambda$. Then the following are equivalent:
\begin{itemize}
\item[(a)] $(\mathcal{U})$ is a maximal ideal of $R$ for every ultrafilter $\mathcal{U}$ in $\mathcal{B}= \prod_{\lambda \in \Lambda} \mathcal{P}(\max(D_\lambda))$.
\item[(b)] The subspace $\max(D_\lambda)$ is proconstructible in $\text{spec}(D_\lambda)$ for every $\lambda \in \Lambda$.
\item[(c)] $D_\lambda$ satisfies property (++) for every $\lambda \in \Lambda$, i.e., for every $r \in D_\lambda$ there exists $d \in D_\lambda$ such that $\mathscr D(r) = \mathscr V(d)$.
\end{itemize}
\end{2.19}

In the particular case where $|\Lambda| = 1$, we get the following statement:

\newtheorem{2.20}[2.1]{Corollary}
\begin{2.20}\label{2.20}
Let $D$ be a commutative ring. Then $\max(D)$ is proconstructible in $\text{spec}(D)$ if and only if $D$ satisfies property (++), i.e.,  for every $r \in D$ there exists $d \in D$ such that $\mathscr D(r) = \mathscr V(d)$.
\end{2.20}

\subsection*{Minimal prime ideals}

For the remainder of this section, let $(D_\lambda)_{\lambda \in \Lambda}$ be a family of integral domains and $R = \prod D_\lambda$. Recall that for an element $x \in R$ we set $z(x) = \{\lambda \in \Lambda \mid x_\lambda = 0 \}$ and define $n(x) = \Lambda \setminus z(x)$. Also, recall the definition of the proper ideal \[(0)_\mathcal{F}^T = \{x \in T \mid z(x) \in \mathcal{F} \}\] depending on a subring $T \subseteq R$ and a filter $\mathcal{F}$ on $\Lambda$.

\theoremstyle{definition}
\newtheorem{2.21}[2.1]{Proposition}
\begin{2.21}\label{2.21}
Let $\mathcal{F}$ be an ultrafilter on $\Lambda$ and $T$ a subring of $R = \prod_{\lambda \in \Lambda} D_\lambda$.
\begin{itemize}
\item[(1)] The ultraproduct $R^* = \prod_{\lambda \in \Lambda}^\mathcal{F} D_\lambda$ is isomorphic to $R/(0)_\mathcal{F}^R$.
\item[(2)] $(0)_\mathcal{F}^T$ is a prime ideal of $T$. 
\item[(3)] Every minimal prime ideal of $T$ is of the form $(0)_\mathcal{F}^T$ for some ultrafilter $\mathcal F$ on $T$.
\end{itemize}
\end{2.21}

\begin{proof}
To prove (1), note that $\varphi: R \to R^*$ mapping an element $r \in R$ to its equivalence class $r^* \in R^*$ is a surjective homomorphism whose kernel is $(0)_\mathcal{F}^R$. 

Now, to prove (2), consider the map $\varphi: T \to R/(0)_\mathcal{F}^R$ with $\varphi(x):= x + (0)_\mathcal{F}^R$. It is a homomorphism whose kernel is $(0)_\mathcal{F}^T$ and whose image is a domain, because, by (1) and the Theorem of \L o\'s, $R/(0)_\mathcal{F}$ is a domain.

Finally, for the proof of (3), let $P \subseteq T$ be a minimal prime ideal and let $M = \{ n(x) \mid x \in T \setminus P \}$. We claim that $M$ has the finite intersection property. Assume to the contrary that there are $x_1,\hdots,x_n \in T \setminus P$ such that $n(x_1) \cap \hdots \cap n(x_n) = \emptyset$. Then $x_1 \cdot \hdots \cdot x_n = 0 \in P$ and therefore there exists some $i \in \{1,\hdots,n\}$ such that $x_i \in P$, which is a contradiction. Let $\mathcal{F}$ be an ultrafilter on $\Lambda$ such that $M \subseteq \mathcal{F}$. Clearly, $T \setminus P \subseteq T \setminus (0)_\mathcal{F}^T$ and therefore $(0)_\mathcal{F}^T \subseteq P$. By the minimality of $P$ it follows that $P = (0)_\mathcal{F}^T$.
\end{proof}

In the next lemma, we restrict our scope to subrings $T \subseteq R$ such that for every $Z \subseteq \Lambda$ there exists some $x \in T$ such that $Z = z(x)$. Note that there are examples of such rings $T$ apart from products of commutative rings, for instance the subring of $R$ generated by all elements $x \in R$ such that $x_\lambda \in \{0,1\}$ for all $\lambda \in \Lambda$.

\newtheorem{2.22}[2.1]{Lemma}
\begin{2.22}\label{2.22}
Let $T$ be a subring of $R = \prod D_\lambda$ such that for every $Z \subseteq \Lambda$ there exists some $x \in T$ such that $Z = z(x)$. Then $(0)_\mathcal{F}^T$ is a minimal prime ideal of $T$ for every ultrafilter $\mathcal{F}$ on $\Lambda$. Moreover, if $\mathcal{F}$ and $\mathcal{G}$ are different ultrafilters on $\Lambda$ then $(0)_\mathcal{F}^T \neq (0)_\mathcal{G}^T$.
\end{2.22}

\begin{proof}
Let $\mathcal{F}$ be an ultrafilter on $\Lambda$ and $P \subseteq T$ be a prime ideal such that $P \subseteq (0)_\mathcal{F}^T$. Let $x \in (0)_\mathcal{F}^T$ and $y \in T$ such that $z(y) = \Lambda \setminus z(x)$. Then $x \cdot y =0 \in P$. Since $P$ is prime, either $x \in P$ or $y \in P$. But $y $ cannot be an element of $P \subseteq (0)_\mathcal{F}^T$, because otherwise $x+y \in (0)_\mathcal{F}^T$ and therefore $\emptyset = z(x+y) \in \mathcal{F}$, which is a contradiction. Hence $x \in P$.

Now, let $\mathcal{G}$ be an ultrafilter on $\Lambda$ different from $\mathcal{F}$. Let $Z \in \mathcal{G} \setminus \mathcal{F}$ and $x \in T$ such that $z(x) = Z$. Then $x \in (0)_\mathcal{G}^T \setminus (0)_\mathcal{F}^T$.
\end{proof}

\newtheorem{2.23}[2.1]{Corollary}
\begin{2.23}\label{2.23}
Let $(D_\lambda)_{\lambda \in \Lambda}$ be a family of integral domains and let $T$ be a subring of $R = \prod_{\lambda \in \Lambda} D_\lambda$ such that for every $Z \subseteq \Lambda$ there exists some $x \in T$ such that $Z = z(x) = \{\lambda \in \Lambda \mid x_\lambda = 0 \}$. Then the map $ \mathcal{F} \mapsto (0)_\mathcal{F}^T$ is a bijection of ultrafilters on $\Lambda$ and minimal prime ideals of $T$.
\end{2.23}

\newtheorem{2.23a}[2.1]{Definition}
\begin{2.23a}\label{2.23a}
Let $\mathcal{U}$ be an ultrafilter in $\mathcal{B}$ and for every $Y \in \mathcal{U}$ set $F_Y = \{\lambda \in \Lambda \mid Y_\lambda \neq \emptyset \}$. We set 
\begin{align*}
\mathcal{F}_\mathcal{U} = \{ F_Y \mid Y \in \mathcal{U} \}
\end{align*}
and call $\mathcal{F}_\mathcal{U}$ the ultrafilter on $\Lambda$ corresponding to $\mathcal{U}$.
\end{2.23a}

Note that $\mathcal{F}_\mathcal{U}$ as defined above is indeed an ultrafilter on $\Lambda$.

\theoremstyle{definition}
\newtheorem{2.24Lemma}[2.1]{Lemma}
\begin{2.24Lemma}\label{2.24Lemma}
Let $R = \prod_{\lambda \in \Lambda} D_\lambda$.
Let $\mathcal{U}$ be an ultrafilter in $\mathcal{B}$ and $\mathcal{F}$ an ultrafilter on $\Lambda$. Then $(0)_\mathcal{F} \subseteq (\mathcal{U})$ if and only if $\mathcal{F} = \mathcal{F}_\mathcal{U}$.
\end{2.24Lemma}

\begin{proof}
Assume that $(0)_\mathcal{F} \subseteq (\mathcal{U})$ and let $F \in \mathcal{F}$. For a subset $M \subseteq \Lambda$ we denote by $\chi_M$ the element of $R$ for which the entry at $\lambda \in \Lambda$ is $1$ if $\lambda \in M$ and is $0$ if $\lambda \notin M$. If we set $M = \Lambda \setminus F$, then $\chi_M \in (0)_\mathcal{F} \subseteq (\mathcal{U})$. Therefore $Y := S(\chi_M) \in \mathcal{U}$, where $Y_\lambda = \emptyset$ if $\lambda \notin F$ and $Y_\lambda = \max(D_\lambda)$ if $\lambda \in F$. So $F = \{ \lambda \in \Lambda \mid Y_\lambda \neq \emptyset \}$ and therefore $F \in \mathcal{F_U}$. Hence $\mathcal{F} \subseteq \mathcal{F_U}$, which implies $\mathcal{F} = \mathcal{F_U}$, because $\mathcal{F}$ is an ultrafilter.

Conversely, let $\mathcal{F} = \mathcal{F_U}$ and let $r \in (0)_\mathcal{F}$. Then for $M = \{\lambda \in \Lambda \mid r_\lambda \neq 0 \}$ we have that $\chi_M \in (0)_\mathcal{F}$. Therefore $\Lambda \setminus M \in \mathcal{F}$, which implies that $\Lambda \setminus M = \{ \lambda \in \Lambda \mid Y_\lambda \neq \emptyset \}$ for some $Y \in \mathcal{U}$. Clearly, we have $S(\chi_M) \geq Y$, so $S(\chi_M) \in \mathcal{U}$. Consequently, $r = r \cdot \chi_M \in (\mathcal{U})$.
\end{proof}

In the following two propositions, we summarize, before we specialize to the case of Prüfer domains, what we can say about the prime spectrum of a product of arbitrary domains.

\theoremstyle{definition}
\newtheorem{2.24}[2.1]{Proposition}
\begin{2.24}\label{2.24}
Let $R = \prod_{\lambda \in \Lambda} D_\lambda$.
\begin{itemize}
\item[(1)] For every prime ideal $\mathfrak{P}$ of $R$ there exists an ultrafilter $\mathcal{U}$ in $\mathcal{B}$ such that $ \mathfrak{P} \subseteq (\mathcal{U})$.
\item[(2)] Every prime ideal $\mathfrak{P}$ of $R$ contains a unique minimal prime ideal. It is equal to $(0)_\mathcal{F_U}$ for every ultrafilter $\mathcal{U}$ in $\mathcal{B}$ such that $\mathfrak{P} \subseteq (\mathcal{U})$.
\end{itemize}
In particular, for every prime $\mathfrak{P}$ or $R$ there exists an ultrafilter $\mathcal{U}$ in $\mathcal{B}$ such that \[(0)_\mathcal{F_U} \subseteq \mathfrak{P} \subseteq (\mathcal{U}).\]

\end{2.24}

\begin{proof}
Let $\mathfrak{P} \subseteq R$ be a prime ideal. Then $\mathfrak{P}$ contains a minimal prime ideal. If $\mathfrak{Q} \subseteq \mathfrak{P}$ is a minimal prime ideal, then by Proposition \ref{2.2} there exists an ultrafilter $\mathcal{F}$ on $\Lambda$ such that $\mathfrak{Q} = (0)_\mathcal{F}$. In the same way, if $\mathfrak{M}$ is a maximal ideal containing $\mathfrak{P}$, then by Proposition \ref{2.21}(3) we can pick some ultrafilter $\mathcal{U}$ on $\mathcal{B}$ such that $\mathfrak{M} = (\mathcal{U})$. Since $(0)_\mathcal{F} \subseteq \mathfrak{P} \subseteq (\mathcal{U})$, it follows that $\mathcal{F} = \mathcal{F_U}$, so $(0)_\mathcal{F}= (0)_\mathcal{F_U}$. Therefore $(0)_\mathcal{F_U}$ is the unique minimal prime ideal contained in $\mathfrak{P}$.
\end{proof}

\newtheorem{2.25}[2.1]{Proposition}
\begin{2.25}\label{2.25}
Let $R = \prod_{\lambda \in \Lambda} D_\lambda$. For an ultrafilter $\mathcal{U}$ in $\mathcal{B}$, let $R^* = \prod_{\lambda \in \Lambda}^\mathcal{F_U} D_\lambda$ be the ultraproduct with respect to $\mathcal{F_U}$ and $R^*_\mathcal{U}$ the localization of $R^*$ at the maximal ideal $(\mathcal{U})^*$ corresponding to~$(\mathcal{U})$.
 
The prime ideals $\mathfrak{P}$ of $R$ with $(0)_\mathcal{F_U} \subseteq \mathfrak{P} \subseteq (\mathcal{U})$ are in inclusion preserving one-to-one correspondence with the prime ideals of $R^*_\mathcal{U}$.

\end{2.25}

\begin{proof}
This follows immediately from the fact that $R^*= R/(0)_\mathcal{F_U}$ and the bijective correspondence between primes of $R^*$ contained in $(\mathcal{U})^*$ and primes of the localization $R^*_\mathcal{U}$.
\end{proof}

The well-known fact that the spectrum of a product $\prod_{\lambda \in \Lambda} K_\lambda$ of fields $K_\lambda$ can be described as the Stone-\v Cech compactification of the discrete space $\Lambda$ is an easy consequence of Corollary~\ref{2.3}, Corollary \ref{2.23} and Proposition \ref{2.24}. Indeed, the Boolean algebra $\mathcal{B} = \prod_{\lambda \in \Lambda} \mathcal{P}(\max(K_\lambda))$ is canonically isomorphic to $\mathcal{P}(\Lambda)$, because each $\max(K_\lambda)$ is a singleton. Thus, in this particular case, ultrafilters in $\mathcal{B}$ are in canonical one-to-one correspondence with ultrafilters on $\Lambda$.

In general, Proposition \ref{2.24} gives some insight into the order structure of the set of prime ideals in the product $R$ of integral domains: $\text{spec}(R)$ is a disjoint union of partially ordered sets $O$, where each $O$ has a unique minimal element. This is also a starting point for our considerations in the next section.

\section{Prime ideals in products of Prüfer domains}\label{section:primeideals}

We now set out to characterize all prime ideals of $R = \prod_{\lambda \in \Lambda} D_\lambda$, where each $D_\lambda$ is a Prüfer domain. By Proposition \ref{2.24}, it is sufficient to describe for every ultrafilter $\mathcal{U}$ in $\mathcal{B}$ the prime ideals $\mathfrak{P} \subseteq R$ with $(0)_\mathcal{F_U} \subseteq \mathfrak{P} \subseteq (\mathcal{U})$. In order to do this, it is more convenient to have a different view upon $\mathcal{B}$.

\newtheorem{3.1}{Remark}[section]
\begin{3.1}\label{3.1Remark}
If $(V_\lambda)_{\lambda \in \Lambda}$ is a family of sets, we denote by $\dot\bigcup_{\lambda \in \Lambda} V_\lambda$ their disjoint union. The map 
\begin{align*}
\mathcal{B} &\to \mathcal{P}\left(\dot\bigcup_{\lambda \in \Lambda} \max(D_\lambda)\right) \\ (V_\lambda)_{\lambda \in \Lambda} &\mapsto \dot{\bigcup}_\lambda V_\lambda 
\end{align*}
defines an isomorphism of Boolean algebras. Using this, we can view ultrafilters in $\mathcal{B}$ as ultrafilters on the set $\dot{\bigcup}_{\lambda \in \Lambda} \max(D_\lambda)$ and the notions defined in Section \ref{section:maxideals} translate as expected. For instance, for an ultrafilter $\mathcal{U}$ on $\dot{\bigcup}_{\lambda \in \Lambda} \max(D_\lambda)$, we define
\[ (\mathcal{U}) = \{ r \in R \mid \{M \in \dot{\bigcup}_{\lambda \in \Lambda} \max(D_\lambda) \mid r_\lambda \in M \} \in \mathcal{U} \}.\]
Analogously, for $U \in \mathcal{U}$, we set $F_U = \{\lambda \in \Lambda \mid U \cap \max(D_\lambda) \neq \emptyset\}$. Then we get an ultrafilter
\[ \mathcal{F_U} = \{ F_U \mid U \in \mathcal{U} \} \]
on $\Lambda$ induced by $\mathcal{U}$. Note that this is the same ultrafilter $\mathcal{F_U}$ that we defined in Section \ref{section:maxideals}, where $\mathcal{U}$ was assumed to be an ultrafilter in $\mathcal{B}$.
\end{3.1}

To characterize all prime ideals of $R$, it suffices, by Proposition~\ref{2.24}, to consider, for every ultrafilter $\mathcal{U}$ on $\dot{\bigcup}_{\lambda \in \Lambda} \max(D_\lambda)$ such that $(\mathcal{U})$ is a maximal ideal, all those primes $\mathfrak{P}$ such that $(0)_\mathcal{F_U} \subseteq \mathfrak{P} \subseteq (\mathcal{U})$.

%
%

From now on, fix an ultrafilter $\mathcal{U}$ on $\dot{\bigcup}_{\lambda \in \Lambda} \max(D_\lambda)$ and let $\mathcal{F_U}$ be the corresponding ultrafilter on $\Lambda$ (as in Remark~\ref{3.1Remark}).
Let $R^* = \prod_{\lambda \in \Lambda}^\mathcal{F_U} D_\lambda$ be the ultraproduct with respect to $\mathcal{F_U}$ and $R^*_\mathcal{U}$ the localization of $R^*$ at the maximal ideal $(\mathcal{U})^*$ corresponding to $(\mathcal{U})$.

For $r \in R$, we denote by $r^*$ the image of $r$ in $R^*$ or the image of $r$ in $R^*_\mathcal{U}$, depending on the context.


\newtheorem{3.1Lemma}[3.1]{Lemma}
\begin{3.1Lemma}\label{3.1Lemma}\cite[Proposition 2.2]{Olb1}
Let $(D_\lambda)_{\lambda \in \Lambda}$ be a family of rings and $\mathcal{F}$ an ultrafilter on $\Lambda$. Then the ultraproduct $R^* = \prod^\mathcal{F}_{\lambda \in \Lambda} D_\lambda$ is a Prüfer domain if and only if \[\{\lambda \in \Lambda \mid D_\lambda \text{ is a Prüfer domain}\} \in \mathcal{F}.\]
\end{3.1Lemma}

\begin{proof}
By the Theorem of \L o\'s, it suffices to express "Prüfer domain" as a sentence in the first-order language of rings. It is clear how to express "domain". For a domain the Prüfer property is equivalent to every non-zero two-generated ideal being invertible, which can be expressed as
\begin{align*}
&\forall a,b \ ((a \neq 0 \lor b \neq 0) \\
 \Rightarrow \ &\exists \alpha, \beta,c,s,t,u,v \
(c \neq 0 \land a \alpha + b \beta = c  \land a\alpha = cs \land b\beta = ct \land a\beta = cu \land b \alpha = cv)),
\end{align*}
see for instance the textbook by Ershov~\cite[Theorem 2.1.1 and Remark 2.1.2]{Ershov}.
\end{proof}

From now on, let $D_\lambda$ be a Prüfer domain for every $\lambda \in \Lambda$. Then $R^*$ is a Prüfer domain and $R^*_\mathcal{U}$ is a valuation domain. Let $K^*$ be the quotient field of $R^*$ and note that it is isomorphic to the ultraproduct $\prod_{\lambda}^\mathcal{F_U} K_\lambda$ of the quotient fields $K_\lambda$ of $D_\lambda$. Moreover we extend the notation $\mathscr V(r_\lambda)$ and $\mathscr D(r_\lambda)$ to elements $k_\lambda \in K_\lambda$, setting \[\mathscr V(k_\lambda) = \{ M \in \max(D_\lambda) \mid v_M(k_\lambda) > 0 \},\] where $v_M$ is the valuation on $K_\lambda$ corresponding to $M$, and $\mathscr D(k_\lambda) = \max(D_\lambda) \setminus \mathscr V(k_\lambda)$. 

In the following, we are able to describe the valuation $v$ on $K^*$ that has $R^*_\mathcal{U}$ as its valuation ring.

\subsection*{Valuations and prime ideals}

For a ring $D$ and a prime ideal $P$ of $D$, as usual, we denote the localization of $D$ at $P$ as $D_P$.

\newtheorem{3.2}[3.1]{Proposition}
\begin{3.2}\label{3.1}
Let $\mathcal{U}$ be an ultrafilter in $\mathcal{B}$.
\begin{itemize}
\item[(1)] The injective homomorphism
\begin{align*}
R^* &\to \prod\nolimits_{M \in \dot\bigcup_{\lambda \in \Lambda} \max(D_\lambda)}^\mathcal{U} (D_\lambda)_M \\
r^* &\mapsto (r_\lambda)_{M \in \dot\bigcup_{\lambda \in \Lambda} \max(D_\lambda)}^\mathcal{U}
\end{align*}
canonically extends to an injective homomorphism 
\[ R^*_\mathcal{U} \to \prod\nolimits^\mathcal{U}_{M \in \dot\bigcup_{\lambda \in \Lambda} \max(D_\lambda)} (D_\lambda)_M\]
via localization.

\item[(2)] For each $M \in \dot\bigcup_{\lambda \in \Lambda} \max(D_\lambda)$, we denote by $\Gamma_M$ the value group of the valuation $v_M$ on $K_\lambda$ with valuation ring $(D_\lambda)_M$.  The map
\begin{align*}
v: K^* \setminus \{0\} &\to \Gamma_\mathcal{U} := \prod\nolimits^\mathcal{U}_{M \in \dot\bigcup_{\lambda \in \Lambda} \max(D_\lambda)} \Gamma_M \\
(k_\lambda)_\lambda^\mathcal{F_U} &\mapsto (v_M(k_\lambda))_M^\mathcal{U}
\end{align*}
is a valuation on $K^*$ with valuation ring $R^*_\mathcal{U}$.

\item[(3)] For $a,b \in R = \prod_\lambda D_\lambda$, the following assertions are equivalent:
\begin{itemize}
\item[(a)] $v(a^*) \geq v(b^*)$.
\item[(b)] $\{M \in \dot\bigcup_{\lambda \in \Lambda} \max(D_\lambda) \mid v_M(a_\lambda) \geq v_M(b_\lambda)\} \in \mathcal{U}$.
\end{itemize}
\end{itemize}

\end{3.2}

\begin{proof}
(1) follows immediately from the definitions of $(\mathcal{U})$ and $\mathcal{F_U}$.

(2) follows from (1) by the simple observation that, in general, an ultraproduct of valuations is a valuation whose valuation ring is the ultraproduct of the corresponding valuation rings. Note that being a valuation ring is a statement in the first-order language of rings.

(3) follows from (2) by definition.
\end{proof}


From now on, we denote by $\Gamma_M$ the value group of the valuation $v_M$ on $K_\lambda$ with valuation ring $(D_\lambda)_M$ and  $\Gamma_\mathcal{U} = \prod^\mathcal{U}_{M \in \dot\bigcup_{\lambda \in \Lambda} \max(D_\lambda)} \Gamma_M$. While $\mathcal{U}$ is fixed, we write $\Gamma$ for $\Gamma_\mathcal{U}$.
We write $S = \{ g \in \Gamma \cup \{\infty\} \mid g > 0\}$. By the definition of elements of ultraproducts as equivalence classes of vectors, we may represent every element of $S$ by a vector $(g_M)$ such that $g_M>0$ for all $M$, thus
\[ S = \{(g_M)_{M \in \dot\bigcup_{\lambda \in \Lambda} \max(D_\lambda)}^{\mathcal{F_U}} \mid \forall M \in \dot\bigcup_{\lambda \in \Lambda} \max(D_\lambda) \ g_M >0 \}.\]

For simplicity, when working with elements $g =(g_M)_{M \in \dot\bigcup_{\lambda \in \Lambda} \max(D_\lambda)}^{\mathcal{F_U}} \in S$, we will therefore always assume that $g_M >0$ for all $M \in  \dot\bigcup_{\lambda \in \Lambda} \max(D_\lambda)$.

Given $g \in S$, we define 
\begin{align*}
(\mathcal{U})^g = \{ x \in R \mid  \exists n \in \mathbb{N} \ v((x^*)^n) \geq g \},
\end{align*}
where $v$ is the valuation as in Proposition~\ref{3.1}.

\newtheorem{3.1a}[3.1]{Proposition}
\begin{3.1a}\label{3.2}
For any $g \in S$, the set $(\mathcal{U})^g$ is a prime ideal of $R$ contained in $(\mathcal{U})$.
\end{3.1a}

\begin{proof}
First, let $g \in S$. Clearly, $(\mathcal{U})^g$ is an ideal. To see that it is contained in $(\mathcal{U})$, let $x \in (\mathcal{U})^g$ and let $n \in \N$ such that $v((x^*)^n) \geq g$. Choose $Y = \{ M \in \dot\bigcup_{\lambda \in \Lambda} \max(D_\lambda) \mid v(x_\lambda^n)\geq g_M \}$. Note that $Y \in \mathcal{U}$ and $g_M > 0$ for all $M$. It follows that $\{M \in \dot\bigcup_{\lambda \in \Lambda} \max(D_\lambda) \mid x_\lambda^n \in M\}  \supseteq Y \in \mathcal{U}$, so $x^n \in (\mathcal{U})$, which is a prime ideal. Therefore $x \in (\mathcal{U})$. \\
Finally, let $a,b \in R$ such that $ab \in (\mathcal{U})^g$ and let  $n \in \mathbb{N}$ such that $v((a^* b^*)^n) \geq g$. Set $Y = \{ M \in \dot\bigcup_{\lambda \in \Lambda} \max(D_\lambda) \mid v_M((a_\lambda b_\lambda)^n) \geq g_M\}$ and note that $Y \in \mathcal{U}$. Given $M \in Y$, it follows that $g_M + g_M \leq v_M(a_\lambda^n b_\lambda^n) + v_M(a_\lambda^n b_\lambda^n) = v_M(a_\lambda^{2n} b_\lambda^{2n}) = v_M(a_\lambda^{2n}) + v_M(b_\lambda^{2n})$. Hence $v_M(a_\lambda^{2n}) \geq g_M$ or $v_M(b_\lambda^{2n}) \geq g_M$. Since $\mathcal{U}$ is an ultrafilter, it follows that the set of those $M$ such that $v_M(a_\lambda^{2n}) \geq g_M$ is in $\mathcal{U}$ or set of those $M$ such that $v_M(b_\lambda^{2n}) \geq g_M$ is in $\mathcal{U}$. Say, the former. Then there exists $n' \in \mathbb{N}$ (namely $n' = 2n$) such that $\{ M \in \dot\bigcup_{\lambda \in \Lambda} \max(D_\lambda) \mid  v_M(a_\lambda^{n'}) \geq g_M \} \in \mathcal{U}$, which means per definition that $a \in (\mathcal{U})^g$.

\end{proof}

\newtheorem{3.3}[3.1]{Proposition}
\begin{3.3}\label{3.3}
Let $x \in (\mathcal{U})$ and $g_M = v_M(x_\lambda)$ if $v_M(x_\lambda) > 0$ and $g_M \in \Gamma_M$ arbitrary otherwise. Then $(\mathcal{U})^g$ is a prime ideal $\mathfrak{P}$ with $(0)_\mathcal{F_U} \subseteq \mathfrak{P} \subseteq (\mathcal{U})$ containing~$x$ and every prime ideal satisfying these conditions contains $(\mathcal{U})^g$.
\end{3.3}

\begin{proof}
We already know that $(\mathcal{U})^g \subseteq (\mathcal{U})$. To see that $x \in (\mathcal{U})^g$, note that $v((x^*)^1) = v(x^*) \geq g$ because $U := \{ M \in \dot\bigcup_{\lambda \in \Lambda} \max(D_\lambda) \mid v_M(x_\lambda) = g_M\} \in \mathcal{U}$ by definition of $(\mathcal{U})$.

Now let $(0)_\mathcal{F_U} \subseteq \mathfrak{P} \subseteq (\mathcal{U})$ be a prime ideal containing $x$. Since the prime ideals of $R$ containing $(0)_\mathcal{F_U}$ and contained in $(\mathcal{U})$ are in inclusion preserving bijective correspondence with the prime ideals of $R^*_\mathcal{U}$, it suffices to prove the inclusion $((\mathcal{U})^g)^* \subseteq \mathfrak{P}^*$ of the corresponding prime ideals in $R^*_\mathcal{U}$. So let $r \in (\mathcal{U})^g$. We show that $r^* \in \mathfrak{P}^*$. Let $n \in \mathbb{N}$ such that $v((r^*)^n) \geq g$ and define $Y$ as the set of those $M$ such that $v_M(r_\lambda^n) \geq g_M$. Note that $Y \in \mathcal{U}$. By Proposition \ref{3.1}, it follows that $v((r^*)^n) \geq v(x^*)$. Therefore, $x^*$ divides $(r^*)^n$ in $R^*_\mathcal{U}$, which implies that $(r^*)^n \in \mathfrak{P}^*$. As $\mathfrak{P}^*$ is a prime ideal, it contains $r^*$.
\end{proof}

\newtheorem{3.4}[3.1]{Theorem}
\begin{3.4} \label{3.4}
Let $R = \prod_{\lambda \in \Lambda} D_\lambda$ where every $D_\lambda$ is a Prüfer domain.

The prime ideals of $R$ are exactly the unions of sets of prime ideals of the form $(\mathcal{U})^g$, where $\mathcal{U}$ is an ultrafilter on $\dot\bigcup_{\lambda \in \Lambda} \max(D_\lambda)$ and $g \in \Gamma_\mathcal{U} \cup \{\infty\}$ with $g > 0$. 


\end{3.4}

\begin{proof}
Since $R^*_\mathcal{U}$ is a valuation domain, the prime ideals of $R$ contained in $(\mathcal{U})$ form a chain. Every union of $(\mathcal{U})^g$ is, therefore, a union of a chain of prime ideals, and hence prime.

Conversely, let $(0)_\mathcal{F_U} \subseteq \mathfrak{P} \subseteq (\mathcal{U})$ be a prime ideal of $R$. For $x \in \mathfrak{P}$ we define an element $g(x) \in S = \{g \in \Gamma_\mathcal{U} \mid g >0\}$. Namely, set $g(x)_M= v_M(x_\lambda)$ if $M \in \mathscr V(x_\lambda)$ and $g(x)_M\in \Gamma_M$ arbitrary otherwise. We claim that 
\begin{align*}
\mathfrak{P} = \bigcup_{x \in \mathfrak{P}} (\mathcal{U})^{g(x)}.
\end{align*}
By Proposition \ref{3.3}, $(\mathcal{U})^{g(x)}$ is the smallest prime ideal contained in $(\mathcal{U})$ and containing~$x$. So $\bigcup_{x \in \mathfrak{P}} (\mathcal{U})^{g(x)} \subseteq \mathfrak{P}$. On the other hand, if $y \in \mathfrak{P}$, then by Proposition \ref{3.3} we have that $y \in (\mathcal{U})^{g(y)}$ and therefore $y \in \bigcup_{x \in \mathfrak{P}} (\mathcal{U})^{g(x)}$. 
\end{proof}

\subsection*{When are the $(\mathcal{U})^g$ distinct?}
In Theorem~\ref{3.4}, we described prime ideals of a product of Prüfer domains as unions of chains of prime ideals of the form $(\mathcal{U})^g$ for a fixed ultrafilter $\mathcal{U}$. It is apparent that many of the $(\mathcal{U})^g$ are superfluous in these unions, because $(\mathcal{U})^g$ and $(\mathcal{U})^h$ can be equal even if $g \neq h$. We now investigate when this is the case by defining an equivalence relation on $S \subseteq \Gamma_\mathcal{U}$.

Moreover, forming the Dedekind-MacNeille completion of the partially ordered set of equivalence classes, which is again a totally ordered set, we get rid of unions of chains of prime ideals and we see that prime ideals correspond to elements of this set. That this correspondence might, in general, still not be bijective is due to the fact that, depending on the component rings $D_\lambda$, the valuation $v$ of Proposition~\ref{3.1} might not be surjective.

\newtheorem{3.1d}[3.1]{Definition}
\begin{3.1d} \label{3.1d}
We now define a relation $\ll$ on $S = \{g \in \Gamma_\mathcal{U} \cup \{\infty\} \mid g >0\}$ (see Proposition~\ref{3.1}), where $\mathcal{U}$ is a fixed ultrafilter on $\dot\bigcup_{\lambda \in \Lambda} \max(D_\lambda)$. Let
\begin{align*}
g\ll h :\Leftrightarrow \forall n \in \N \ n \cdot g < h
\end{align*}
for $g,h \in S$. 
An equivalence relation on $S$ is given by
\[ g \sim h :\Leftrightarrow \exists n \in \N \ (n \cdot g \geq h \land n \cdot h \geq g), \]
for $g,h \in S$. We write $[g]_\sim$ for the equivalence class of $g \in S$.
\end{3.1d}

Note that $g \ll h$ depends only on $[g]_\sim$ and $[h]_\sim$ and, therefore, $\ll$ also defines a total order on the set $S/_\sim$ of equivalence classes of $\sim$. 

\theoremstyle{definition}
\newtheorem{3.5}[3.1]{Lemma}
\begin{3.5} \label{3.5}
Let $g,h \in S$.
\begin{itemize}
\item[(1)] If $(\mathcal{U})^h \subsetneqq (\mathcal{U})^g$, then $g \ll h$.
\item[(2)] Assume that for every $(w_M)_{M \in \dot\bigcup_{\lambda \in \Lambda} \max(D_\lambda)} \in \prod \Gamma_M$ there exists $x \in R = \prod D_\lambda$ such that $\{M \in \dot\bigcup_{\lambda \in \Lambda} \max(D_\lambda) \mid v_M(x_\lambda) = w_M\} \in \mathcal{U}$. 
Then $g\ll h$ implies $(\mathcal{U})^h \subsetneqq (\mathcal{U})^g$.
\end{itemize}
Note that, in particular, (2) is the case if the set of all $\lambda$ such that $D_\lambda$ is semilocal is in $\mathcal{F_U}$.
\end{3.5}

\begin{proof}
To see (1), let $x \in (\mathcal{U})^g \setminus (\mathcal{U})^h$. Then there exists  $n \in \mathbb{N}$ such that $v((x^*)^n) \geq g$. On the other hand, $v((x^*)^{n'}) < h$ for all $n' \in \N$. Let $m \in \N$. Then $m \cdot g \leq v((x^*)^{nm}) < h$. It follows that $g \ll h$. 

Now let $g \ll h$.
By the additional assumption in statement (2), we can find $x \in R$ with $(v_M(x_\lambda))_M^{\mathcal{U}} = g$. Then $x \in (\mathcal{U})^g \setminus (\mathcal{U})^h$.

\end{proof}

Lemma~\ref{3.5}(1), shows that the definition of $(\mathcal{U})^g$ is independent from the choice of $g$ as long as $g$ lies in the same equivalence class with respect to $\sim$. Thus, we define 
\[(\mathcal{U})^\mathfrak{g} = (\mathcal{U})^g\]
for $\mathfrak{g} = [g]_\sim \in S/_\sim$. Now, $(\mathcal{U})^\mathfrak{g} \neq (\mathcal{U)}^\mathfrak{h}$ if and only if $\mathfrak{g} \neq \mathfrak{h}$, under the assumption of Lemma~\ref{3.5}(2).

\newtheorem{3.1r}[3.1]{Remark}
\begin{3.1r}\label{3.1r}

Instead of unions of prime ideals of the form $(\mathcal{U})^\mathfrak{g}$, we can use elements of a Dedekind--MacNeille completion of $S/_\sim$ to describe all prime ideals of a product of Prüfer domains. Let $\Delta$ be any partially ordered set. A \textit{Dedekind--MacNeille completion} of $\Delta$ is a complete lattice $\overline{\Delta}$ together with an order embedding $ \iota: \Delta \to \overline{\Delta}$ such that every order embedding of $\Delta$ into a complete lattice factors uniquely through $\iota$.

Every element of a Dedekind--MacNeille completion of $\Delta$ is an infimum (resp. a supremum) of some subset of $\Delta$.
Every partially ordered set $\Delta$ has a Dedekind--MacNeille completion, given by the set of all Dedekind cuts of $\Delta$. If $\Delta$ is totally ordered then so is every Dedekind--MacNeille completion.
\end{3.1r}

We fix a Dedekind--MacNeille completion $\overline{S/_\sim}$ of the totally ordered set $S/_\sim$. For $\mathfrak{g} \in \overline{S/_\sim}$, we define
\begin{align*}
(\mathcal{U})^\mathfrak{g} = \{ x \in R \mid  \exists n \in \mathbb{N} \ [v((x^*)^n)]_\sim \geq \mathfrak{g} \}.
\end{align*}
Note that $\mathfrak{g} = \inf \{[g]_\sim \mid g \in S \land [g]_\sim \geq \mathfrak{g}\}$ and therefore, as a union of a chain of prime ideals,
\[(\mathcal{U})^\mathfrak{g} = \bigcup_{\substack{g \in S \\ [g]_\sim \geq \mathfrak{g}}} (\mathcal{U})^g\]
is a prime ideal. Using Theorem~\ref{3.4}, this immediately leads to the following alternative description of the prime spectrum of a product of Prüfer domains.

\newtheorem{Corollary}[3.1]{Corollary}
\begin{Corollary}\label{Corollary}
Let $R = \prod_{\lambda \in \Lambda} D_\lambda$ where every $D_\lambda$ is a Prüfer domain.


The prime ideals of $R$ are exactly the sets of the form $(\mathcal{U})^\mathfrak{g}$, where $\mathcal{U}$ is an ultrafilter on $\dot\bigcup_{\lambda \in \Lambda} \max(D_\lambda)$ and $\mathfrak{g} \in \overline{S/_\sim}$ where $S =  \{g \in \Gamma_\mathcal{U} \cup \{\infty\} \mid g >0\}$, $\sim$ is as in Definition~\ref{3.1d} and $\overline{S/_\sim}$ is the Dedekind--MacNeill completion as in Remark~\ref{3.1r}.
\end{Corollary}

\subsection*{Heights of prime ideals}
We now introduce a special type of ultrafilter that will be helpful to force certain prime ideals in $R$ to have infinite height. 
An ultrafilter $\mathcal{F}$ on a set $\Lambda$ is \textit{countably incomplete} if there exists a countable partition $(P_n)_{n \in \mathbb{N}}$ of $\Lambda$ such that $P_n \notin \mathcal{F}$ for every $n \in \mathbb{N}$. It is shown in~\cite[Theorem 6.1.4]{mod} that for every infinite set $\Lambda$ there exists a countably incomplete ultrafilter on $\Lambda$. 

\newtheorem{3.6}[3.1]{Lemma}
\begin{3.6} \label{3.6}
If $\mathcal{F}$ is a countably incomplete ultrafilter on $\Lambda$, then every ultrafilter $\mathcal{V}$ on $\dot\bigcup_{\lambda \in \Lambda} \max(D_\lambda)$ with $\mathcal{F} = \mathcal{F_U}$ is countably incomplete.
\end{3.6}

\begin{proof}
Let $(P_n)_{n \in \mathbb{N}}$ be a partition of $\Lambda$ such that $P_n \notin \mathcal{F}$ for all $n \in \mathbb{N}$. Define $Q^{(n)} \subseteq \dot\bigcup_{\lambda \in \Lambda} \max(D_\lambda)$ for each $n \in \mathbb{N}$ such that $Q^{(n)} \cap \max(D_\lambda) = \max(D_\lambda)$ if $\lambda \in P_n$ and $Q^{(n)} \cap \max(D_\lambda) = \emptyset$ else. Clearly, $\bigcup_{n \in \mathbb{N}}Q^{(n)} = \dot\bigcup_{\lambda \in \Lambda} \max(D_\lambda)$ and $Q^{(m)} \cap Q^{(n)} = \emptyset$ for all $m, n \in \mathbb{N}$ with $m \neq n$.

Assume that $Q^{(n)} \in \mathcal{V}$ for some $n \in \mathbb{N}$. Then $P_n = \{\lambda \in \Lambda \mid Q^{(n)}\cap \max(D_\lambda) \neq \emptyset \} \in \mathcal{F_V} = \mathcal{F}$, which is a contradiction.
\end{proof}

Recall that we fixed an ultrafilter $\mathcal{U}$ on $\dot\bigcup_{\lambda \in \Lambda} \max(D_\lambda)$. Let $\mathcal{F_U}$ the induced ultrafilter on $\Lambda$.

\newtheorem{3.7}[3.1]{Lemma}
\begin{3.7} \label{3.7}
Assume that $\mathcal{F_U}$ is countably incomplete and let $g,h \in \Gamma \cup \{\infty\}$ with $g, h >0$.
\begin{itemize}
\item[(1)] If $g \ll h$, then there exists some $k \in S$ such that $g \ll k \ll h$.
\item[(2)] If $g < \infty$, then there exists some $k \in S$ such that $g \ll k < \infty$. 
\end{itemize}
\end{3.7}

\begin{proof}
(2) follows immediately from (1) by setting $h = \infty$. Note that $g \ll \infty$ if and only if $g < \infty$.

To show (1), we can assume without loss of generality that $g_M < h_M$ for all $M \in \dot\bigcup_{\lambda \in \Lambda} \max(D_\lambda)$.
We define two complementary subsets $V,W$ of $\dot\bigcup_{\lambda \in \Lambda} \max(D_\lambda)$ as follows:
\begin{align*}
V \cap \max(D_\lambda) &= \{M \in \max(D_\lambda) \mid \forall n \in \mathbb{N} \ n \cdot g_M < h_M\} \\
W \cap \max(D_\lambda) &= \{ M \in \max(D_\lambda) \mid \exists N \in \mathbb{N} \ N \cdot g_M \geq h_M\}.
\end{align*}
Since $\mathcal{U}$ is an ultrafilter on $\dot\bigcup_{\lambda \in \Lambda} \max(D_\lambda)$, we either have $V \in \mathcal{U}$ or $W \in \mathcal{U}$.
Assume that $V \in \mathcal{U}$. Since $\mathcal{U}$ is countably incomplete by Lemma \ref{3.6}, we can choose a partition $(P_n)_{n \in \mathbb{N}}$ of $\dot\bigcup_{\lambda \in \Lambda} \max(D_\lambda)$ such that $P_n \notin \mathcal{U}$ for all $  n \in \mathbb{N}$. By setting $V^{(n)} = P_n \cap V$ for each $n \in \mathbb{N}$, we get a partition $(V^{(n)})_{n \in \mathbb{N}}$ of $V$ such that $V^{(n)} \notin \mathcal{U}$ for every $n \in \mathbb{N}$. For $M \in \dot\bigcup_{\lambda \in \Lambda} \max(D_\lambda)$, we define $k_M = n \cdot g_M$ if $M \in V_\lambda^{(n)}$ and $k_M = g_M$ if $M \notin V_\lambda$. Then clearly $k \ll h$.

To see that $g \ll k$, let $n \in \mathbb{N}$. Note that $n \cdot g < k$ is equivalent to $ \lnot (n \cdot g \geq k)$, which is the case if and only if $\forall U \in \mathcal{U} \ \exists M \in U \ n \cdot g_M < k_M$.

To show this last statement, let $U \in \mathcal{U}$. Then $U \cap V \in \mathcal{U}$ and therefore there exists some $N>n$ such that $U \cap V^{(N)} \neq \emptyset$. (Otherwise $U \cap V^{(N)} = \emptyset$ for all $N>n$ and, therefore, $(U \cap V^{(1)}) \cup \hdots \cup (U \cap V^{(n)}) = U \cap V \in \mathcal{U}$, which implies $U \cap V^{(i)} \in \mathcal{U}$ for some $ i$, a contradiction.) Pick some $M \in U \cap V^{(N)}$. Then $k_M = N \cdot g_M > n \cdot g_M$. This shows that $g \ll k$.

Now consider the case where $W \in \mathcal{U}$. We define $k$ as follows. For each $M \in W$, there exists some $N > 1$ such that $h_M \leq N \cdot g_M$ and we pick $N_M \geq 1$ such that $N_M \cdot g_M < h_M \leq (N_M + 1) \cdot g_M$ and define $k_M = [N_M/\log(N_M)] \cdot g_M$, where $[\, . \, ]$ denotes the floor function and $[\infty] := \infty$. For $M \notin W$, we define $k_M = g_M$.

Let $n \in \mathbb{N}$. Using equivalent statements to $n \cdot g < k$ and $n \cdot k< h$ as above, we see that it suffices to show, for every $Y \in \mathcal{U}$:
\begin{itemize}
\item[(1)] There exists $M \in Y$ such that $n \cdot g_M < k_M$.
\item[(2)] There exists $M \in Y$ such that $n \cdot k_M < h_M$.
\end{itemize}
So let $Y \in \mathcal U$. We can assume without loss of generality that $Y \subseteq W$. First, note that the set $\{N_M \mid M \in Y \}$ is unbounded, because $g \ll h$. It follows that the sets
\begin{align*}
S_{1} &:= \{ [N_M/\log(N_M)] \mid M \in Y\} \text{ and} \\
S_{2} &:= \{ \log(N_M) \mid M \in Y\}
\end{align*}
are also unbounded.
To show (1), we use that $S_{1}$ is unbounded and pick $M \in Y$ such that $n < [N_M/\log(N_M)]$. It follows that $n \cdot g_M < [N_M/\log(N_M)] \cdot g_M = k_M$. \\
For the proof of (2), we can pick $M \in Y$ such that $n < \log(N_M)$. It follows that $n \cdot [N_M/\log(N_M)] \leq n \cdot N_M/\log(N_M) < n \cdot N_M/n = N_M$. Hence $n \cdot k_M = n \cdot [N_M/\log(N_M)] \cdot g_M < N_M \cdot g_M < h_M$.
\end{proof}

Note that the assumptions of the following theorem are always satisfied if all $D_\lambda$ are semilocal Prüfer domains with a uniform bound on the cardinalities of $\max(D_\lambda)$).

\newtheorem{3.9}[3.1]{Theorem}
\begin{3.9} \label{3.9}
Let $(D_\lambda)_{\lambda \in \Lambda}$ be a family of Prüfer domains and set $R = \prod_{\lambda \in \Lambda} D_\lambda$. Assume that $\mathcal{F}$ is a countably incomplete ultrafilter on $\Lambda$ and $\mathcal{U}$ an ultrafilter on $\dot\bigcup_{\lambda \in \Lambda} \max(D_\lambda)$ such that $\mathcal{F}$ is the induced ultrafilter on $\Lambda$. 

If $\mathcal{U}$ has the additional property that there exist $Y \in \mathcal{U}$ and $N \in \mathbb{N}$ with $|Y \cap \max(D_\lambda)| \leq N$ for all $\lambda \in \Lambda$ then for every prime ideal $P \subseteq R$ with $(0)_\mathcal{F} \subsetneqq P \subseteq (\mathcal{U})$ there exists some prime ideal $Q \subseteq R$ such that $(0)_\mathcal{F} \subsetneqq Q \subsetneqq P$.

In particular, every prime ideal $P$ of $R$ with $(0)_\mathcal{F} \subsetneqq P \subseteq (\mathcal{U})$ is of infinite height.
\end{3.9}

\begin{proof}
Let $I \subseteq S$ such that $P = \bigcup_{g \in I} (\mathcal{U})^g$, which exists by Theorem \ref{3.4}. Since $(0)_\mathcal{F} \subseteq (\mathcal{U})^g$ for all $g \in I$, there must exist some $g \in I$ such that $(\mathcal{U})^\infty = (0)_\mathcal{F} \subsetneqq (\mathcal{U})^g$. It follows by Lemma~\ref{3.5}(1) that $g < \infty$. So by Lemma~\ref{3.7}(2) there exists some $h \in S$ with $g \ll h < \infty$. Lemma~\ref{3.5}(2) implies that $(0)_\mathcal{F} \subsetneqq (\mathcal{U})^h \subsetneqq (\mathcal{U})^g \subseteq P$. The assertion follows by setting $Q = (\mathcal{U})^h$.
\end{proof}


\section*{Acknowledgements}
We would like to thank the anonymous referee for the very helpful comments on our paper; they helped to improve the exposition crucially.

\nocite{Heinzer1}
\nocite{Heinzer2}
\nocite{Heinzer3}

\bibliographystyle{amsplainurl}
\bibliography{bibliography}

\vspace{0.5cm}
\noindent
\textsc{Carmelo A. Finocchiaro, Dipartimento
di Matematica e Informatica, Università degli Studi di Catania, 95125 Catania, Italy} \\
\textit{E-mail address}: \texttt{cafinocchiaro@unict.it}\\

\noindent
\textsc{Sophie Frisch, Department of Analysis and Number Theory (5010), Technische Universität Graz, Kopernikusgasse 24, 8010 Graz, Austria} \\
\textit{E-mail address}: \texttt{frisch@math.tugraz.at} \\

\noindent
\textsc{Daniel Windisch, Department of Analysis and Number Theory (5010), Technische Universität Graz, Kopernikusgasse 24, 8010 Graz, Austria} \\
\textit{E-mail address}: \texttt{dwindisch@math.tugraz.at}

\end{document}